\documentclass[smallextended,numbook,runningheads]{svjour3}     
\pdfoutput=1

\smartqed  
\usepackage{graphicx}
\usepackage{amssymb, amsmath, latexsym, amsfonts, mathcomp}
\usepackage{color}
\usepackage{xcolor}
\usepackage{mathptmx}      
\usepackage{xfrac}
\usepackage{algorithm}
\usepackage{algpseudocode}

\usepackage[backend=bibtex]{biblatex}
\def\S2{{S^2}}

\def\dd{\mathrm{d}}

\newcommand{\R}{\mathbb{R}}

\newcommand{\g}{\mathfrak{g}}

\journalname{BIT}
\begin{document}

\title{Variational symplectic diagonally implicit Runge-Kutta methods for isospectral systems\thanks{Funded by the Deutsche Forschungsgemeinschaft (DFG, German Research Foundation) – Project-ID 422037413 – TRR 287. 
}}


\author{Clauson Carvalho da Silva         \and
        Christian Lessig 
}


\institute{C. C. da Silva, C. Lessig \at
           Otto-von-Guericke-Universit{\"a}t Magdeburg, Institut f{\"u}r Simulation und Graphik \\
              \email{\url{clauson@isg.cs.uni-magdeburg.de}, \url{christian.lessig@ovgu.de}}           
}

\date{Received: date / Accepted: date}

\maketitle

\begin{abstract}
Isospectral flows appear in a variety of applications, e.g. the Toda lattice in solid state physics or in discrete models for two-dimensional hydrodynamics, with the isospectral property often corresponding to mathematically or physically important conservation laws.
Their most prominent feature, i.e. the conservation of the eigenvalues of the matrix state variable, should therefore be retained when discretizing these systems.
Recently, it was shown how isospectral Runge-Kutta methods can, in the Lie-Poisson case also considered in our work, be obtained through Hamiltonian reduction of symplectic Runge-Kutta methods on the cotangent bundle of a Lie group.
We provide the Lagrangian analogue and, in the case of symplectic diagonal implicit Runge-Kutta methods, derive the methods through a discrete Euler-Poincar{\'e} reduction.
Our derivation relies on a formulation of diagonally implicit isospectral Runge-Kutta methods in terms of the Cayley transform, generalizing earlier work that showed this for the implicit midpoint rule.
Our work is also a generalization of earlier variational Lie group integrators that, interestingly, appear when these are interpreted as update equations for intermediate time points.
From a practical point of view, our results allow for a simple implementation of higher order isospectral methods and we demonstrate this with numerical experiments where both the isospectral property and energy are conserved to high accuracy. 

\keywords{Runge Kutta methods \and Lie group integrators \and variational integrators}
\subclass{MSC 65L06 \and MSC 65P10}
\end{abstract}

\section{Introduction}
\label{sec:intro}

Isospectral systems, for which the time evolution preserves the spectrum of a matrix state variable, arise in the mathematical description of various systems of interest.
They are generically given by
\begin{align}
  \label{eq:isospectral_conceputal}
  \dot{A} = \big[ B(A) , A \big] , \quad A(0) = A_0
\end{align}
where $A$, $B(A)$ are matrices and the spectrum of $A$ is preserved~\cite[IV.3.2]{hairer2006geometric}.
See~\cite{calvo1997numerical} for a detailed discussion of the properties of these systems.

A classical example for an isospectral system is the ideal rigid body in its reduced formulation~\cite{Marsden1999}, as described by Euler's equations. 
In this case the angular momentum (as an anti-symmetric matrix in $\mathfrak{so}^*(3)$) is the isospectrally conserved quantity. 
The periodic Toda lattice in its Lax pair formulation is also isospectral with the isospectral state $L$ in this case being symmetric and the generator of the dynamics $B(L)$ being anti-symmetric.
Another example is the discrete 2D barotropic vorticity equation on the torus and 2-sphere when the ``matrix model'', i.e. the representation theory of the discrete Heisenberg group respectively $\mathrm{SO}(3)$, is used for spatial discretization~\cite{Zeitlin1991,Hoppe1998,Zeitlin2004}. 
The isospectral property is in this case a discrete analogue (with $N$ conserved quantities for $N$ degrees of freedom) of the integrated powers of vorticity that are conserved in the continuous system~\cite{Khesin1998}.

Since conserved quantities are an important intrinsic characteristic of a physical system, it is natural to preserve isospectrality under discretization. 
An early work that accomplished this was those by Moser and Veselov~\cite{Moser1991} for rigid body like systems. 
Moore, Mahony, and Helmke~\cite{Moore1994} introduced isospectral time integrators for the solution of problems in linear algebra, since many classical ones there, such as the QR or eigen decomposition, can be cast as isospectral flows.
Isospectral numerical time integrators were also considered by Diele, Lopez and Politi who proposed the use of the Cayley transform~\cite{Diele1998}.
Calvo, Iserles and Zanna~\cite{calvo1997numerical} showed that Runge-Kutta methods applied to Eq.~\ref{eq:isospectral_conceputal} break isospectrality.
To overcome this, they introduce modified Gauss-Legendre Runge-Kutta methods that provide isospectral time integration schemes of arbitrary high order.
Recently, Modin and Vivani~\cite{modin2019lie} developed an alternative route to overcome the obstruction from~\cite{calvo1997numerical} by applying a Runge-Kutta scheme to the cotangent bundle $T^*G$ of a quadratic matrix Lie group $G$ and then reducing to the dual Lie algebra $\g^*$ using the momentum map, analogous to how Lie-Poisson reduction works in the continuous case, cf.~\cite[Ch. 13]{Marsden1999}.
The ansatz exploits that many isospectral flows are Lie-Poisson, i.e. their phase space is a dual Lie algebra $\mathfrak{g}^*$, and the dynamics are reduced ones from $T^*G$.
Vivani~\cite{viviani2019minimal} recently extended~\cite{modin2019lie} and showed that the implicit midpoint (Runge-Kutta) scheme can be written and implemented efficiently with an intermediate time step (in the spirit of the Verlet method).
He also pointed out that the Cayley transform then arises naturally and does not have to be introduced as a group map as in~\cite{Diele1998}.

Another line of research on discretizations that retain conserved quantities is based on the Lagrangian formulation of mechanics and uses a discrete Hamilton's principle~\cite{marsden2001discrete}.
This ansatz has also been developed for reduced dynamical systems on Lie groups and algebras~\cite{BouRabee2007,BouRabee2009}.
Gawlik et al.~\cite{gawlik2011geometric}, in particular, developed a formulation for semi-direct products that, as a special case, is applicable to systems such as the rigid body.

In the present work, we will derive the Lagrangian analogue of the work of Modin and Vivani~\cite{modin2019lie} and show that isospectral Runge-Kutta methods can be obtained from the reduced Hamilton's principle of Bou-Rabee and Marsden~\cite{BouRabee2009} and Gawlik et al.~\cite{gawlik2011geometric}.
We do so for the case of symplectic diagonally implicit Runge-Kutta methods (SDIRKs), which are an important special case since they are the simplest and most easily implemented symplectic Runge-Kutta schemes.
Similar to~\cite{modin2019lie}, we will consider isospectral Lie-Poisson systems for our derivation.
We also extend the results of Vivani~\cite{viviani2019minimal} by showing that the Cayley transform arises for an SDIRK applied at the group level when it is written using intermediate time steps located between the Runge-Kutta steps.
Another result of the present work is that the midpoint rule by Gawlik et al.~\cite{gawlik2011geometric} contains, in disguise, an isospectral method; more precisely, the scheme from~\cite{gawlik2011geometric} updates the midpoints of the isospectral midpoint rule of Vivani~\cite{viviani2019minimal}.
We exemplify our results by implementing our isospectral SDIRKs for the rigid body, the Toda lattice, and the matrix model discretization of the 2D Euler fluid.
We show that with a $7$-stage $4^{\textrm{th}}$-order SDIRK, the oscillations in the energy, which are standard for a symplectic integrator, are on the order of machine precision. 
Within the accuracy afforded by standard double precision, the integrator therefore preserves the invariants of the continuous system.
SDIRKs provide here the advantage that these can be implemented more easily and can be computed more efficiently than general isospectral Runge-Kutta methods.

The remained of the paper is structured as follows.
In Sec.~\ref{sec:preliminaries} we will recall the necessary background on Lie groups, their algebras and isospectral Lie-Poisson systems.
In Sec.~\ref{sec:rk_cayley} we will show that the Cayley transform arises for any SDIRK on a quadratic matrix Lie group by introducing intermediate time steps.
The central result of the paper, namely a variational derivation of isospectral SDIRKs, 
will be presented in Sec.~\ref{sec:variational}.
Numerical results for the rigid body, the Toda lattice, and the discrete 2D Euler fluid are discussed in Sec.~\ref{sec:numerics}.

%

\section{Preliminaries}
\label{sec:preliminaries}

\subsection{Lie groups, Lie algebras and retraction maps}
\label{sec:preliminaries:lie}

Let $G$ be a matrix Lie group with identity $e$ and $\mathfrak{g}$ its Lie algebra.
We will denote the conjugate transpose of an element $g \in G$ with $g^{\dagger}$.
If conjugation is defined as a left action $g \cdot h = g \,h \,g^{-1}$, the adjoint action $\mathrm{Ad}_{g} : G \times \g \to \g$ of the group $G$ on its Lie algebra $\g$ is given by $\text{Ad}_{g} \cdot \xi = g\, \xi\, g^{-1}$.
Infinitesimally it becomes $\mathrm{ad}_{\xi}(\eta) = [\xi,\eta]$. 
The dual of $\mathrm{ad}$ is the infinitesimal coadjoint action $\text{ad}^* : \g \times \g^* \to \g^*$ defined through $\smash{\langle \mathrm{ad}_{\xi}( \eta ) , \alpha \rangle = \langle \eta , \mathrm{ad}_{\xi}^* ( \alpha ) \rangle}$, which it plays an important role for Lie-Poisson dynamics.
Unless mentioned otherwise, we will work with the pairing between $\g$ and $\g^*$ associated with the Frobenius inner-product on $\mathfrak{g}$ given by $\langle \xi , \alpha \rangle = \langle \xi , \alpha^\sharp \rangle = \mathrm{tr}( \xi^{\dagger}  \alpha^\sharp )$ which also allows us to identify $\g$ and $\g^*$. We ommit the sharp symbol $\sharp$ throughout the paper.
A local diffeomorphism $\tau:\mathfrak{g} \to G$ satisfying $\tau(0)= e$, $\tau(\xi)^{-1} = \tau(-\xi)$ and whose tangent map fulfills $T_{0}\tau = \text{Id}$ is called a retraction map (or local group diffeomorphism). Such a map provides a local chart for the Lie group around the identity (with coordinates of the first kind) and induces an atlas for $G$ through the translation action. 
The tangent $T_\xi \tau$ of $\tau$ maps $T_\xi \mathfrak{g} \cong T_e G = \mathfrak{g}$ to $T_{\tau(\xi)}G$.
The same is accomplished with the tangent map $T_e R_{\tau(\xi)}$ where $R_{g}$ is right translation $R_g(h) = h \cdot g$. 
The maps $T_\xi \mathfrak{g}$ and $T_e R_{\tau(\xi)}$ are related by the right trivialization $\dd \tau_{\xi}: \mathfrak{g} \to \mathfrak{g}$ of $\tau$,
\begin{align*}
    T_\xi\tau= T_e R_{\tau(\xi)}\circ \dd \tau_{\xi}.     
\end{align*}
Since $\tau$ is a local diffeomorphism, the map $T_\xi \tau$ is invertible for $\xi$ close enough to $0 \in \mathfrak{g}$ which implies that $\dd \tau_{\xi}$ is also invertible.
We will write $\dd\tau^{-1}_{\xi} \equiv (\dd\tau_{\xi})^{-1}: \mathfrak{g} \to \mathfrak{g}$. A basic property of the right trivialization used throughout the paper is $\dd\tau_{\xi} =  \text{Ad}_{\tau(\xi)} \circ \dd\tau^{-1}_{-\xi}$, cf.~\cite{Owren2016}.

In the following, we will consider quadratic matrix Lie groups, i.e. we will work with $G\subset \mathrm{GL}(N)$ so that there exists a $J$ such that
$J^\dagger \, g \, J = g, \forall g$.
For $\xi \in \g$ it then holds that $J\,\xi + \xi^\dagger J = 0$. 
Examples of quadratic Lie groups are $\mathrm{O}(n)$ and $\mathrm{SO}(n)$ with $J = \mathrm{Id}$ and the symplectic group $\mathrm{Sp}(2N,\R)$ with $J$ the symplectic matrix.
An expedient retraction map for quadratic Lie groups is the Cayley transform~\cite{BouRabee2007}
\begin{align}
  \label{eq:cayley}
    \mathrm{cay}(\xi) = \left( \text{Id}-\frac{\xi}{2} \right)^{-1} \! \left(\text{Id}+\frac{\xi}{2} \right) .
\end{align}
The two factors in the definition above commute. 
The right trivialization $\dd \tau$ for the Cayley transform is given by~\cite{hairer2006geometric}
\begin{align*}
    \dd\mathrm{cay}_{\xi}\cdot \eta &= \left( \text{Id}-\frac{\xi}{2} \right)^{-1} \eta \left(\text{Id}+\frac{\xi}{2} \right)^{-1}
    \\[3pt]
    \dd\mathrm{cay}^{-1}_{\xi}\cdot \eta &= \left( \text{Id}-\frac{\xi}{2} \right) \eta \left(\text{Id}+\frac{\xi}{2} \right).
\end{align*}
The configuration space considered in the paper are quadratic Lie groups. 
Therefore the Cayley transform plays a major role throughout the paper.



\subsection{Lie-Poisson systems on reductive Lie algebras}
\label{sec:preliminaries:lie_poisson}

We recall Lie-Poisson systems on the dual of a Lie algebra $\mathfrak{g} \subseteq \mathfrak{gl}(n)$ and the associated reconstruction process to the cotangent bundle of the Lie group $T^*G$. The flow of a Lie-Poisson system on a dual Lie algebra $\mathfrak{g}^*$ is a solution to
\begin{align}\label{eq:lie_poisson}
  \dot{\mu} 
  = \text{ad}^*_{\nabla H(\mu)} \mu, \quad \mu\in\mathfrak{g}^*
\end{align}
where $H : \mathfrak{g}^* \to \R$ is a Hamiltonian function. 
Using the Frobenius inner product 
$\langle \xi, \eta \rangle = \mathrm{tr}(\xi^{\dagger} \eta)$ in $\mathfrak{g}$,
we can identify $\mathfrak{g}$ with $\mathfrak{g}^*$ and the dual of the infinitesimal adjoint action then becomes
$\mathrm{ad}^*_\xi (\eta) = \mathrm{\Pi}([\xi^\dagger, \eta])$
where $\mathrm{\Pi}$ is the projection from $\mathfrak{gl}(n)$ onto $\mathfrak{g}$. For the case where $\mathfrak{g}$ is a reductive Lie algebra, it is known that
$[\mathfrak{g}^{\dagger}, \mathfrak{g}] \subset \mathfrak{g}$ and therefore the Lie-Poisson equation~\ref{eq:lie_poisson} becomes
\begin{align}\label{eq:lie_poisson_algebra}
  \dot{\mu} = \left[ \nabla H(\mu)^{\dagger}, \mu\right].
\end{align}
An important property of such systems is isospectrality, i.e, the eigenvalues of the matrices $\mu(t)$ are independent of $t$. 

Every Lie-Poisson system can be reconstructed to a canonical Hamiltonian system on $T^*G$. The main tool in Lie-Poisson reconstruction is the momentum map
$\mathcal{J}(g,p) : T^*G \to \g^*$ which for the right action  $(g,p)\cdot h = (gh, p(h^{-1})^\dagger)$ of the group $G \subseteq \mathrm{GL}(n)$ on the phase space $T^*G$, is given by $\mu = \mathcal{J}(g,p) = g^{\dagger} p$.
With $\mathcal{J}$ we can define a left invariant Hamiltonian $\tilde{H}(g,p)= H(\mathcal{J}(g,p)) = H(g^{\dagger} p)$ and the canonical Hamilton equations on $T^*G$ are
\begin{align}
  \label{eq:hamilton_eq_iso}
  \begin{aligned}
  \dot{g} &= g \nabla H(g^\dagger p) = g B(g^\dagger p)^\dagger  
  \\[4pt] 
  \dot{p} &= -p \nabla H(g^\dagger p)^\dagger = - p B(g^\dagger p)
  \end{aligned}
\end{align}
where we write $B(\cdot) = \nabla H (\cdot)^\dagger: \mathfrak{g} \to \mathfrak{g}$ as in \cite{modin2019lie}. For the right invariant version of the Lie-Poisson reconstruction, see Appendix~\ref{sec:right_invariance}.
\subsection{Symplectic diagonally implicit Runge-Kutta methods}
\label{sec:preliminaries:sdirk}

The Butcher tableau for an $s$-stage symplectic diagonally implicit Runge-Kutta methods (SDIRK) is~\cite[Ch.~V.2.1]{hairer2006geometric}
\begin{align}
  \label{eq:sdirk:butcher}
\begin{array}{c|cccc}
 c_1 & \frac{b_1}{2}  &   \\[3pt]
 c_2 & b_1  & \frac{b_2}{2} &  &  \\[3pt]
 \vdots &  \vdots &  \vdots & \ddots \\[3pt]
 c_s & b_1  & b_2  & \cdots & \frac{b_s}{2}\\[3pt]
\hline\rule{0pt}{1.01\normalbaselineskip}
    &  b_1  &  b_2 & \cdots & b_s \\[3pt]
\end{array}
\end{align}
Given a fixed time step $h$, we will write $h_i = h \, b_i$ and $r_i = \sum_{j=1}^i b_j$ for $i =0, \ldots, s$ so that, through the diagonal entries in Eq.~\ref{eq:sdirk:butcher}, $c_i$ and $r_{i-1}$ always differ by $b_i / 2$.
Hence $\sum^s_i h_i = h$ and $r_s = 1$ and also $b_0=0$. 
An important result for us is that an $s$-stage SDIRK can be decomposed in $s$ implicit midpoint rule steps with suitable step size~\cite[Theorem VI.4.4]{hairer2006geometric}.


\section{The Cayley transform and symplectic diagonally implicit Runge-Kutta methods}
\label{sec:rk_cayley}

As discussed in the last section, any $s$-stage symplectic diagonally implicit Runge-Kutta method (SDIRK) can be decomposed into $s$ implicit midpoint rule steps~\cite[Theorem VI.4.4]{hairer2006geometric}.
In this section we will show that this carries over to isospectral SDIRKs and hence the simplicity of the isospectral minimal midpoint algorithm presented in~\cite{viviani2019minimal} also holds for higher order isospectral SDIRKs.
We will also show that the Cayley transform appears naturally in this context.
The key to this result is to introduce not only the usual intermediate time steps that appear in a Runge-Kutta method, i.e. $k_{n,i}^g$ and $k_{n,i}^p$ in the Hamiltonian setting, but also intermediate half time steps at $t_i^n = h n + h \sum_{j=1} b_j$.
These intermediate half points, denoted $g_{n+r_i}$ and $p_{n+r_i}$, can be seen as a generalization of the intermediate step that one considers for the Verlet integrator and they are defined through the strictly lower diagonal part of the Butcher tableau in Eq.~\ref{eq:sdirk:butcher}.

Although the Lie groups we consider are nonlinear, since they are quadratic it nonetheless holds that a symplectic Runge-Kutta method applied linearly to Hamilton's equations in Eq.~\ref{eq:hamilton_eq_iso} yields $g_n \in G, \forall n$.
Exploiting this for an $s$-stage SDIRK, the Runge-Kutta steps can be written in terms of the intermediate Runge-Kutta points as
\begin{subequations}
\label{eq:cayley:der}
\begin{align}
  \label{eq:cayley:der:0}
    k_{n,i}^g 
    &= \left( g_{n,r_{i-1}} + \frac{h_i}{2}  \, k_{n,i}^g \right) B({\mu}_{n,c_i})^\dagger
    \\[4pt]
    \label{eq:cayley:der:01}
    k_{n,i}^p
    &= -\left( p_{n,r_{i-1}} + \frac{h_i}{2}  \, k_{n,i}^p \right) B({\mu}_{n,c_i}) .
\end{align}
where the reduced momentum is 
\begin{align}
  {\mu}_{n,c_i} = \left( g_{n,r_{i-1}} + \frac{h_i}{2}  \, k_{n,i}^g \right)^{\! \dagger} \! \left( p_{n,r_{i-1}} + \frac{h_i}{2}  \, k_{n,i}^p \right)
\end{align}
\end{subequations}
by the definition of the momentum map.
The intermediate points $g_{n+r_i}$ and $p_{n+r_i}$ are defined as
\begin{align*}
  g_{n,r_i} = g_n + \sum_{j=0}^i h_j k_{n,j}^g \quad \quad 
  p_{n,r_i} = p_n + \sum_{j=0}^i h_j k_{n,j}^p 
\end{align*}
and they are always a half step away from $k_{n,i+1}^g$ and $k_{n,i+1}^p$, cf. Eq.~\ref{eq:sdirk:butcher}. Here we have $g_{n,r_0} = g_n$ and $g_{n,r_s} = g_{n+1}$ (and analogous for the variable $p$).  

The first intermediate point at $t = h n + b_1$ is therefore
\begin{align*}
    g_{n,r_1} = g_n + h_1 \, k_{n,1}^g = g_n + h_1 \, \left(g_n + \frac{h_1}{2} k_{n,1}^g \right) B({\mu}_{n,c_1})^\dagger .
\end{align*}
Fully expanding the right hand side of the last equation and using from the left hand side that $k_1^g = ( g_{n,r_1} - g_n) / h_1$ we obtain
\begin{align*}
   g_{n,r_1} \left( \text{Id} - \frac{h_1}{2}B({\mu}_{n,c_1})^{\dagger} \right)  = g_n \left( \text{Id} + \frac{h_1}{2}B({\mu}_{n,c_1})^{\dagger} \right) .
\end{align*}
The last equation states conceptually that taking a (linear) half step backward from $g_{n,r_1}$ and a (linear) half step forward from $g_n$ agree. 
With the definition of the Cayley transform in Eq.~\ref{eq:cayley} and taking into account that the terms in the parenthesis commute, this is equivalently given by
\begin{align*}
    g_{n,r_1} = g_n \, \mathrm{cay}(h_1 \, B({\mu}_{n,c_1})^{\dagger}).
\end{align*}
Following the same steps as above for $g_{n,r_i}$ one obtains 
    $g_{n,r_i} = g_{n,{r_{i-1}}} \, \mathrm{cay}(h_i \, B({\mu}_{n,c_i})^\dagger)$.
and we conclude that
\begin{align}
  \label{eq:cayley:g_n:final}
    g_{n+1} = g_n \prod_{i=1}^s \mathrm{cay}(h_i \, B({\mu}_{n,c_i})^\dagger) .
\end{align}
An analogous calculation shows that 
\begin{align}
  \label{eq:cayley:p_n:final}
    p_{n+1} = p_n \prod_{i=1}^s \mathrm{cay}(h_i \, B({\mu}_{n,c_i}))^{-1} .
\end{align}
With the momentum map $\mathcal{J}$, the phase space points $(q_n,p_n)$ are related to a sequence of reduced momenta $\mu_n \in \g^*$ by $\mu_n = \mathcal{J}(g_n,p_n) = g_n^{\dagger} p_n$.
Using Eq.~\ref{eq:cayley:g_n:final} and Eq.~\ref{eq:cayley:p_n:final} we thus have
\begin{align}
\label{eq:isospectral_property}
    \mu_{n+1} 
    = \left(\prod_{i=s}^1\mathrm{cay}(h_i \, B({\mu}_{n,c_i}))\right)\, \mu_n \left(\prod_{i=1}^s\mathrm{cay}(h_i \, B({\mu}_{n,c_i}))^{-1}\right) .
\end{align}
Using the momentum map to define the reduced momentum variable also at the intermediate point we have
\begin{align} \label{eq:update_intermediate_half}
  \mu_{n,r_i} = g_{n,r_i}^\dagger p_{n,r_i} = \mathrm{cay}(h_{i}\, B({\mu}_{n,c_i})) \mu_{n, r_{i-1}} \mathrm{cay}(h_{i}\, B({\mu}_{n,c_i}))^{-1} .
\end{align}
and a calculation shows that 
\begin{align*}
  \mu_{n,c_i} = \dd \mathrm{cay}_{h_{i}\, B({\mu}_{n,c_i}))}\, \mu_{n, r_{i-1}}
\end{align*}
Eq.~\ref{eq:update_intermediate_half} can be written as
\begin{align*}
   \mu_{n,r_i} = \mathrm{Ad}_{\mathrm{cay}(h_{i}\, B({\mu}_{n,c_i}))} \, \mu_{n, r_{i-1}} = \dd \mathrm{cay}^{-1}_{-h_{i}\, B({\mu}_{n,c_i}))}\,  \circ \dd \mathrm{cay}_{h_{i}\, B({\mu}_{n,c_i}))} \mu_{n, r_{i-1}} .
\end{align*} 
One thus has the following algorithm for updating the intermediate algebra points
\begin{subequations} 
\label{eq:algorithm}
\begin{align}
  \mu_{n, r_{i-1}} &=  \dd \mathrm{cay}^{-1}_{h_{i}\, B({\mu}_{n,c_i}))}\,\mu_{n,c_i}\\
  \mu_{n, r_{i}} &=  \dd \mathrm{cay}^{-1}_{-h_{i}\, B({\mu}_{n,c_i}))}\,\mu_{n,c_i} 
\end{align}
\end{subequations}
where one solves implicitly the first equation to find $\mu_{n,c_i}$ from the known $\mu_{n, r_{i-1}}$ and uses the solution in the second equation to update to $\mu_{n, r_{i}}$. 
The $\mu_{n, r_{i}}$ and $\mu_{n,c_i}$ thus form a staggered grid that can be used in conjunction to advance in time to the full steps.
The existence of solution of the implicit step is proven for sufficiently small $h$ in~\cite[Lemma 3]{viviani2019minimal}. See appendix~\ref{sec:right_invariance} for the right invariant version of the above algorithm.
We summarize this section with the following theorem

\begin{theorem}
  Consider the reduced isospectral Hamiltonian system in Eq.~\ref{eq:lie_poisson_algebra} 
  and the Butcher tableau for a symplectic diagonally implicit Runge-Kutta method in Eq.~\ref{eq:sdirk:butcher}. 
  Then Eq.~\ref{eq:algorithm} provides an isospectral integrator for the system.
\end{theorem}


The theorem follows directly from Eq.~\ref{eq:isospectral_property} and this provides an alternative proof of Theorem~4 in the recent work by Modin and Vivani~\cite{modin2019lie} for the case of a symplectic diagonally implicit Runge-Kutta method.
Previously, Vivani~\cite{viviani2019minimal} already obtained an analogous result for the case of the implicit midpoint rule.
Our extension to SDIRKs partially verifies a conjecture that can be found in this work~\cite[Remark~4]{viviani2019minimal}. Furthermore, we showed that~\cite[Theorem VI.4.4]{hairer2006geometric}, which states that an $s$-stage SDIRK can be decomposed in $s$ implicit midpoint rule steps with suitable step size, also holds in the isospectral reduced case.


\begin{remark}
  \label{rem:averaging_cotangent_bundle}
  In the work by Calvo, Iserles, and Zanna~\cite{calvo1997numerical}, one has that $\mu_{n,1/2} = (\mu_n + \mu_{n+1})/2$, i.e, the midpoint used in the update equation is the average of two consecutive full algebra points. A similar interpretation can be given for the isospectral integrators given by Eq.~\ref{eq:algorithm}.
  Defining the variables
  \begin{subequations}
    \begin{align*}
  	  g_{n,c_i} &= g_{n,r_i} + \frac{h_i}{2} k_{n,i}^g = \frac{g_{n,r_i}+g_{n,r_{i+1}}}{2}\\
  	  p_{n,c_i} &= g_{n, r_i} + \frac{h_i}{2} k_{n,i}^p = \frac{p_{n,r_i}+p_{n,r_{i+1}}}{2}
  	\end{align*}
  \end{subequations}
  we have that the intermediate half points $\mu_{n,c_i} = g_{n,c_i}^\dagger p_{n,c_i}$ result from the averaging of the intermediate points $(g_{n,r_i}, p_{n,r_i})$ and $(g_{n,r_{i+1}}, p_{n,r_{i+1}})$ at the cotangent bundle level.
\end{remark}

\section{Variational symplectic diagonally implicit Runge-Kutta methods}
\label{sec:variational}

From the general theory of reduction one expects that the isospectral symplectic Runge-Kutta methods of Modin and Vivani that were derived through Hamiltonian reduction and the momentum map can also be obtained through a reduced discrete variational principle.
In the following theorem we show that this is indeed the case for isospectral SDIRKs. 
The Cayley transform, which serendipitously appeared in the previous section, will be an integral part for this.
Our result also has a surprising connection to an earlier variational Lie group integrator in the work by Gawlik et al.~\cite{gawlik2011geometric}, see Remark~\ref{remark:connection:gawlik}.

\begin{theorem}
\label{th:variational_generalSDIRK}
Let $G$ be a matrix Lie Group with Lie algebra $\mathfrak{g}$ and $L:TG \to \mathbb{R}$ be a left invariant Lagrangian with left trivialization $\ell:\mathfrak{g} \to \mathbb{R}$ given by
\begin{align*}
  \ell(\xi) = \frac{1}{2} \left \langle B^{-1}(\xi^\dagger), \xi \right \rangle
\end{align*}
where $B:\mathfrak{g}^* \to \mathfrak{g}$ is smooth and invertible.
Let $\{ g_{n,r_1}, g_{n,r_2}, \ldots, g_{n,r_{s}} \}_{n=0,\ldots,N-1} \subset G$ for a positive integer $s$ with $g_n = g_{n,r_0} = g_{n-1,r_s}$
and let $\{\xi_{n,c_1}, \ldots, \xi_{n,c_s}\}_{n=0,\ldots,N-1} \subset \mathfrak{g}$ be the associated  sequence in the Lie algebra such that
\begin{align}
    \label{eq:group_update_generalSDIRK}
    g_{n,r_i} = g_{n,r_{i-1}} \tau \big( h_i \, \xi_{n,c_i} \big)
\end{align}
with $\{b_i, c_i\}$ being the coefficients of the Butcher tableau of an $s$-stage SDIRK scheme.
If the sequence $\{g_{n,r_i}\}$ is stationary for the action $s_d: G^{N} \to \mathbb{R}$
\begin{align*} 
    s_d(\{g_{n,r_i}\}) 
    = \sum_{n=0}^{N-1} \sum_{i = 1}^s h_i \, l( \xi_{n,c_i}) 
\end{align*}
under discrete variations of $\{ g_{n,r_i} \}$ with fixed endpoints, then the sequece $\{\xi_{n,c_i}\}$ and associated $\mu_{n,c_i} = B^{-1}(\xi_{n,c_i}^\dagger) \in \mathfrak{g}^*$ have to satisfy the discrete Euler-Poincar{\'e} equations
\begin{align}
\label{eq:discrete_euler_poincare}
\begin{aligned}
    \left(\dd \tau^{-1}_{h_{i+1} \xi_{n,c_{i+1}}}\right)^* \mu_{n,c_{i+1}} &= \left(\dd \tau^{-1}_{-h_{i} \xi_{n,c_i}}\right)^* \mu_{n,c_i} \qquad i = 1, \ldots, s-1
    \\[3pt]
    \left(\dd \tau^{-1}_{h_1 \,\xi_{n+1,c_1}}\right)^* \mu_{n+1,c_1} &= \left(\dd \tau^{-1}_{-h_s\, \xi_{n,c_s}}\right)^* \mu_{n,c_s} \qquad i = s .
\end{aligned}
\end{align}
Furthermore, in the case of a quadratic Lie group and given a initial condition $\mu_0$, the sequence $\{\mu_n\}$ defined by
\begin{align}
  \label{eq:rk:variational:full_points}
    \mu_{n+1} = \left(\mathrm{Ad}_{\prod_{i=1}^s \tau(h_j \, \xi_{n,c_j})}\right)^*  \mu_n = \mathrm{Ad}^*_{\left(\prod_{i=1}^s \tau(h_j \, \xi_{n,c_j})\right)^{-1}} \, \mu_n.
\end{align}
agrees with the isospectral SDIRK method of the previous section.
\end{theorem}

Computationally, Eq.~\ref{eq:discrete_euler_poincare} allows to implicitly update the intermediate reduced momenta $\mu_{n,c_{i}}$ and with Eq.~\ref{eq:rk:variational:full_points} the next full time step can be obtained.
Eq.~\ref{eq:discrete_euler_poincare} is of Euler-Poincar{\'e} type by considering $\dd \tau_{-h_{i} \xi_{n,c_i}}^{-1}$ as the finite time analogue of the infinitesimal coadjoint action $\mathrm{ad}_{\xi_{n,c_i}}^*$ in the continuous Euler-Poincar{\'e} equations.

\begin{proof}
The stationarity of the discrete action under variations of the group sequence $\{ g_{n,r_i}\}$  with fixed endpoints means
\begin{align*}
    0 = \frac{d}{d\epsilon} \Big|_{\epsilon = 0} s_d\left(\{ g_{n,r_i}^\epsilon \}\right)  
    =\sum_{n=0}^{N-1} \sum_{i = 1}^s\frac{d}{d\epsilon} h_i \, l( \xi^\epsilon_{n,c_i}) \, \Big|_{\epsilon = 0} 
    = \sum_{n=0}^{N-1} \sum_{i = 1}^s h_i \, \left \langle \frac{\delta l}{\delta \xi_{n,c_i}}, \delta \xi_{n,c_i} \right \rangle 
\end{align*}
with $g^\epsilon_{0} = g_0$, $g^\epsilon_{N} = g_N$ and $\{g^0_{n,r_i}\} = \{g_{n,r_i}\}$ for all $n,i, \epsilon$.
 Through Eq.~\ref{eq:group_update_generalSDIRK}, a discrete variation in the group $\{ \delta g_{n,r_i} \} = \{ \frac{d}{d\epsilon}  g^\epsilon_{n,r_i} \big\vert_{\epsilon = 0} \}$ induces a variation in the corresponding algebra sequence. 
Writing $\eta = g^{-1} \delta g $ and using the properties of retraction maps, the induced algebra variations are given by
\begin{align*}
    \delta \xi_{n,c_1} = \frac{1}{h_1} \left( -\dd \tau^{-1}_{h_1 \xi_{n,c_1}}(\eta_{n,r_0}) + \dd \tau^{-1}_{-h_1 \xi_{n,c_1}}(\eta_{n,r_1}) \right).
\end{align*}
Stationarity of the discrete action thus amounts to
\begin{align*}
  0 &= \sum_{n=0}^{N-1} \sum_{i = 1}^s\, \left \langle \frac{\delta l}{\delta \xi_{n,c_i}}, 
    \left( -\dd \tau^{-1}_{h_i \xi_{n,c_i}}(\eta_{n,r_{i-1}}) + \dd \tau^{-1}_{-h_i \xi_{n,c_i}}(\eta_{n,r_{i}}) \right)
    \right\rangle.
\end{align*}
Since $g_0^\epsilon = g_0$ and $g_N^\epsilon = g_N$ for every $\epsilon$, we have $\delta g_0 = \delta g_N = 0$ and $\eta_0 = \eta_N =0$. 
Eq.~\ref{eq:discrete_euler_poincare} then follows by duality and using a standard summation by parts argument. 
From $\mu_{n,c_i} = B^{-1}(\xi_{n,c_i}^\dagger)$ it follows that Eq.~\ref{eq:rk:variational:full_points} and Eq.~\ref{eq:isospectral_property} are equivalent.

\end{proof}

\begin{figure}
  \centering
  \includegraphics[width=\textwidth]{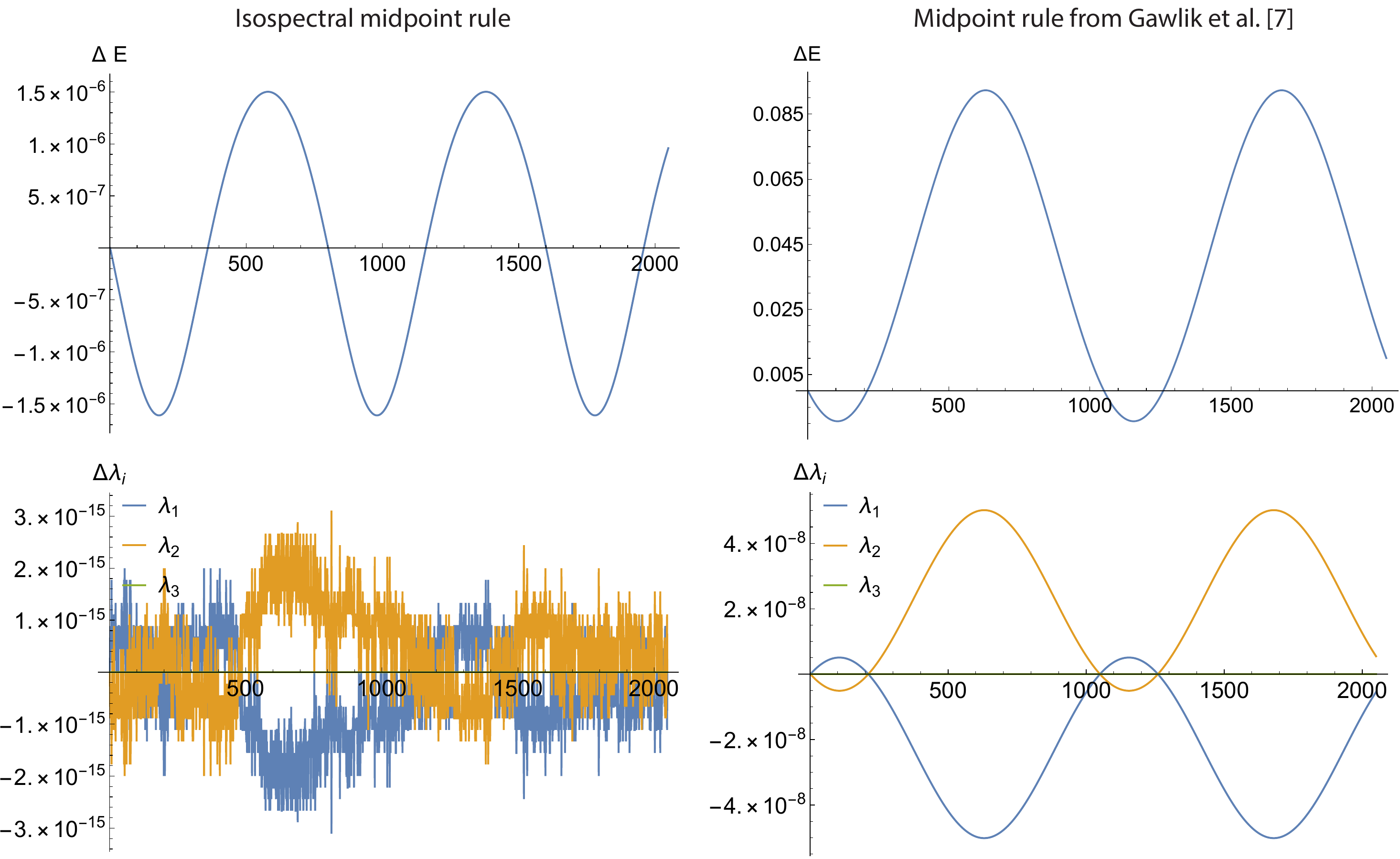}  
  \caption{Comparison of isospectral midpoint rule (left) and the implicit midpoint rule from Gawlik et al.~\cite{gawlik2011geometric} (right) for change in energy (top) and change in the eigenvalues (bottom).}
  \label{fig:gawlik}
\end{figure}

\begin{remark}
  In the continous case, identifying $\mathfrak{g}$ and $\mathfrak{g}^*$, the Euler-Poincar{\'e} equation for an isospectral flow is
  \begin{align*}
  	\dot{\mu} = [\xi, \mu]
  \end{align*}
  The solution $\mu(t)$ and the group curve satisfying $g^{-1}(t) \dot{g}(t) = \xi(t)$ are related by~\cite[IV.3.2]{hairer2006geometric}
  \begin{align*}
  	g(t) \, \mu(t) g(t)^{-1} = g(0) \, \mu(0) g(0)^{-1} \quad \text{ for all $t$}. 
  \end{align*}  
  This ensures that the $\mu(t)$ are related by a similarity transformation for all $t$ and therefore all possess the same spectrum. Eq.~\ref{eq:rk:variational:full_points} are equivalent to the discrete equation
  \begin{align*}
  	g_{n+1} \, \mu_{n+1} g_{n+1}^{-1} = g_n \, \mu_n g_n^{-1} 
  \end{align*}  
  which provides a direct verification that our integrator preserves the isospectral property.  
\end{remark}

\begin{remark}
  \label{remark:connection:gawlik}
  The simplest integrator that results from the above theorem is the implicit midpoint rule studied by Vivani~\cite{viviani2019minimal}.
  Interestingly, it also appeared in disguise in the work by Gawlik et al.~\cite{gawlik2011geometric} on variational discretizations of (possibly infinite dimensional) systems on Lie groups. 
  Combining a forward half-step of the implicit midpoint rule from the intermediate point $\tilde{\mu}_k$ with a backward half-step from $\tilde{\mu}_{k+1}$ one obtains the following implicit update rule for the half-steps
  \begin{align}
    \label{eq:gawlik:update}
    \dd \mathrm{cay}^{-1}_{h B(\tilde{\mu}_{k+1})} \tilde{\mu}_{k+1} &= \dd \mathrm{cay}^{-1}_{-h B(\tilde{\mu}_{k})} \tilde{\mu}_{k} 
  \end{align} 
  or, equivalently,  
  \begin{align*}
   \left( \mathrm{Id} - \frac{h B(\tilde{\mu}_{k+1})}{2} \right) \tilde{\mu}_{k+1} \left( \mathrm{Id} + \frac{h B(\tilde{\mu}_{k+1})}{2}\right) &= \left( \mathrm{Id} + \frac{h B(\tilde{\mu}_{k})}{2} \right) \tilde{\mu}_{k} \left( \mathrm{Id} - \frac{h B(\tilde{\mu}_{k})}{2}\right) .
  \end{align*}
  This are Eq.~4.9 and Eq.~4.12 in the work by Gawlik et al.~\cite{gawlik2011geometric}. 
  However, since only the half time steps are considered in this work their integrator is not isospectral and also yields significantly stronger oscillations in the energy, see Fig.~\ref{fig:gawlik}.    
\end{remark}



\section{Numerical Experiments}
\label{sec:numerics}

In this section we present numerical results for the isospectral SDIRKs from Theorem~\ref{th:variational_generalSDIRK} for the rigid body, the Toda lattice, and the matrix model discretization of the barotropic vorticity equation on $S^2$.
We begin with some details on the computations.
Our implementation is available online at [will be made available upon publication of the manuscript].

\subsection{Implementation of isospectral SDIRKs for quadratic Lie groups}

For the case of isospectral systems defined on quadratic Lie algebras, Eq.~\ref{eq:discrete_euler_poincare} and Eq.~\ref{eq:rk:variational:full_points} are equivalent to Eq.~\ref{eq:algorithm}. 
Given an SDIRK and using the same notation as in Sec.~\ref{sec:preliminaries:sdirk}, Theorem~\ref{th:variational_generalSDIRK} amounts to
\begin{align}
\label{eq:algorithm_matrix}
  \begin{aligned}
  \mu_{n,r_{i-1}} &= \left( \text{Id}-\frac{h_i}{2}B\left( \mu_{n,c_i}\right) \right) \mu_{n,c_i} \left(\text{Id}+\frac{h_i}{2} B\left( \mu_{n,c_i}\right) \right)
 	\\[3pt]
\mu_{n,r_i} &= \left( \text{Id}+\frac{h_i}{2}B\left( \mu_{n,c_i}\right) \right) \mu_{n,c_i} \left(\text{Id}-\frac{h_i}{2} B\left( \mu_{n,c_i}\right) \right)
\end{aligned}
\end{align}   
which allows to update the $\mu_{n,r_i}$ and $\mu_{n,c_i}$ by leap frogging between them.
We summarize the overall computations in Algo.~\ref{alg:isospectral_sdirk}.


In the experiments, we will consider the $2^{\textrm{nd}}$-order implicit midpoint rule and the isospectral integrators obtained from a $2$-stage $2^{\textrm{nd}}$-order SDIRK~\cite{kalogiratou2012diagonally} and a $7$-stage $4^{\textrm{th}}$-order SDIRK~\cite{kalogiratou2012diagonally}.    
   
\begin{figure}
  \includegraphics[width=\textwidth]{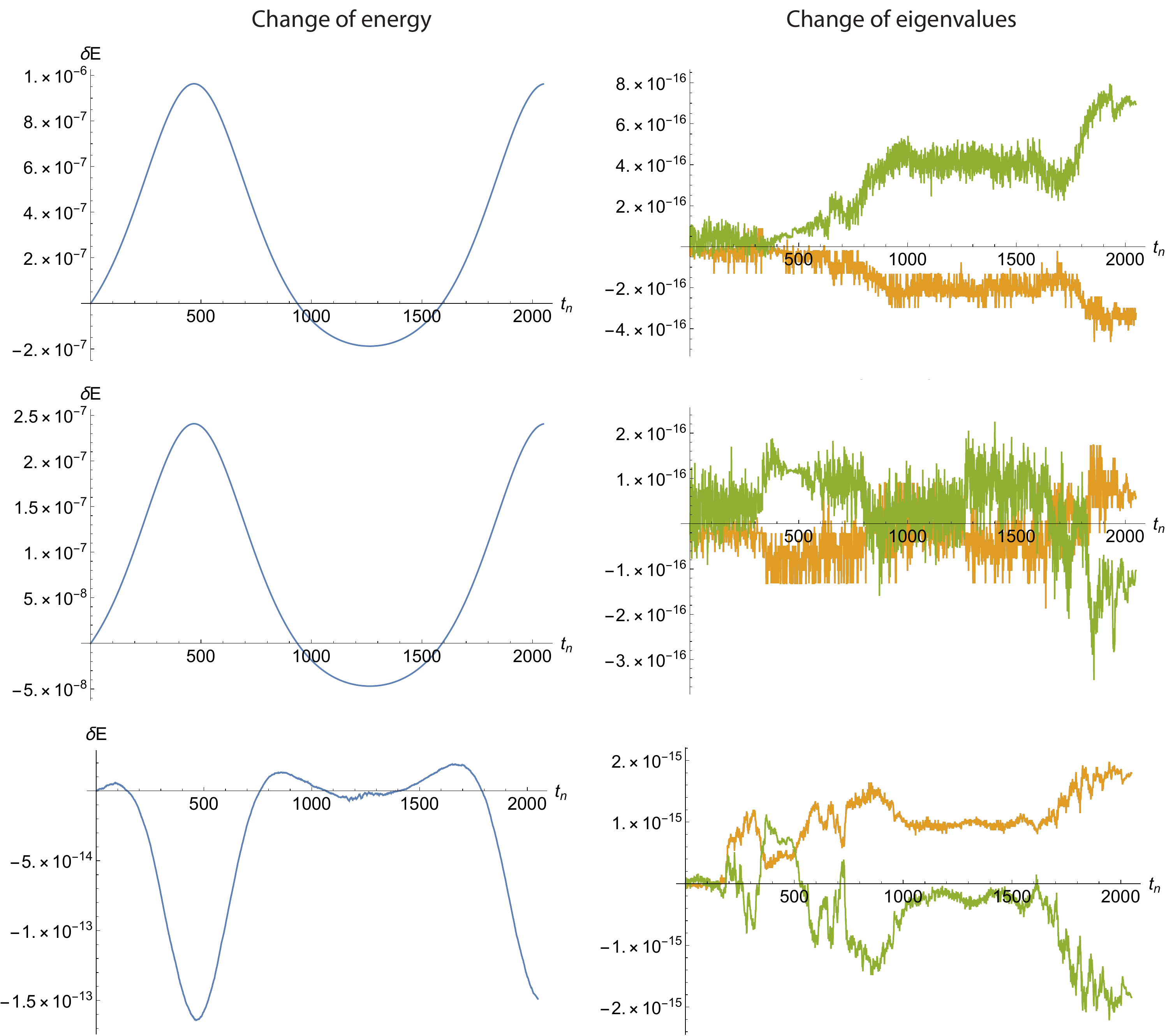}
  \caption{Change in energy and change in eigenvalues for the rigid body for implicit midpoint rule (top), a $2$-stage $2^{\textrm{nd}}$-order scheme (middle) and a $7$-stage $4^{\textrm{nd}}$-order scheme (bottom).}
  \label{fig:rigid_body_experiments}
\end{figure}

\subsection{Rigid body}

The ideal rigid body is the classical example of an isospectral system. 
Its dynamics are described by Euler's equations, which, in matrix form, are~\cite[Ch. 15]{Marsden1999}
\begin{align}
  \dot{W}= [ \nabla H(W)^\dagger, W]
\end{align}
where $W \in \mathfrak{so}^*(3)$ and the Hamiltonian is $H(W) = \frac{1}{2} \left \langle \mathcal{I}^{-1}W, W \right \rangle$ with $\mathcal{I}$ a symmetric $3 \times 3$ matrix.

Fig.~\ref{fig:rigid_body_experiments} shows results for the three isospectral integrators we consider for a random skew-symmetric matrix with entries in $[-1,1]$ as initial condition and with a time step of $h=0.01$.
The superiority of the $4^{\textrm{th}}$-order integrator is evident in the conservation of the Hamiltonian, which is almost on the order of machine precision. 
The conservation of the angular momentum (not shown) is of the same order as those for the eigenvalues.

\subsection{Toda lattice}

The Toda lattice describes the pairwise interactions of particles on a line exerting exponential forces on each other~\cite{Toda1967}. 
When the line is periodic, the Lax pair formulation of the Toda lattice~\cite[Ch. IV.3.2]{hairer2006geometric} results in the isospectral flow $\dot{L} = \left[ B(L), L\right]$ with Hamiltonian $H(L) = 2 \text{Tr}(L^2)$ for 
\begin{figure}
  \includegraphics[width=\textwidth]{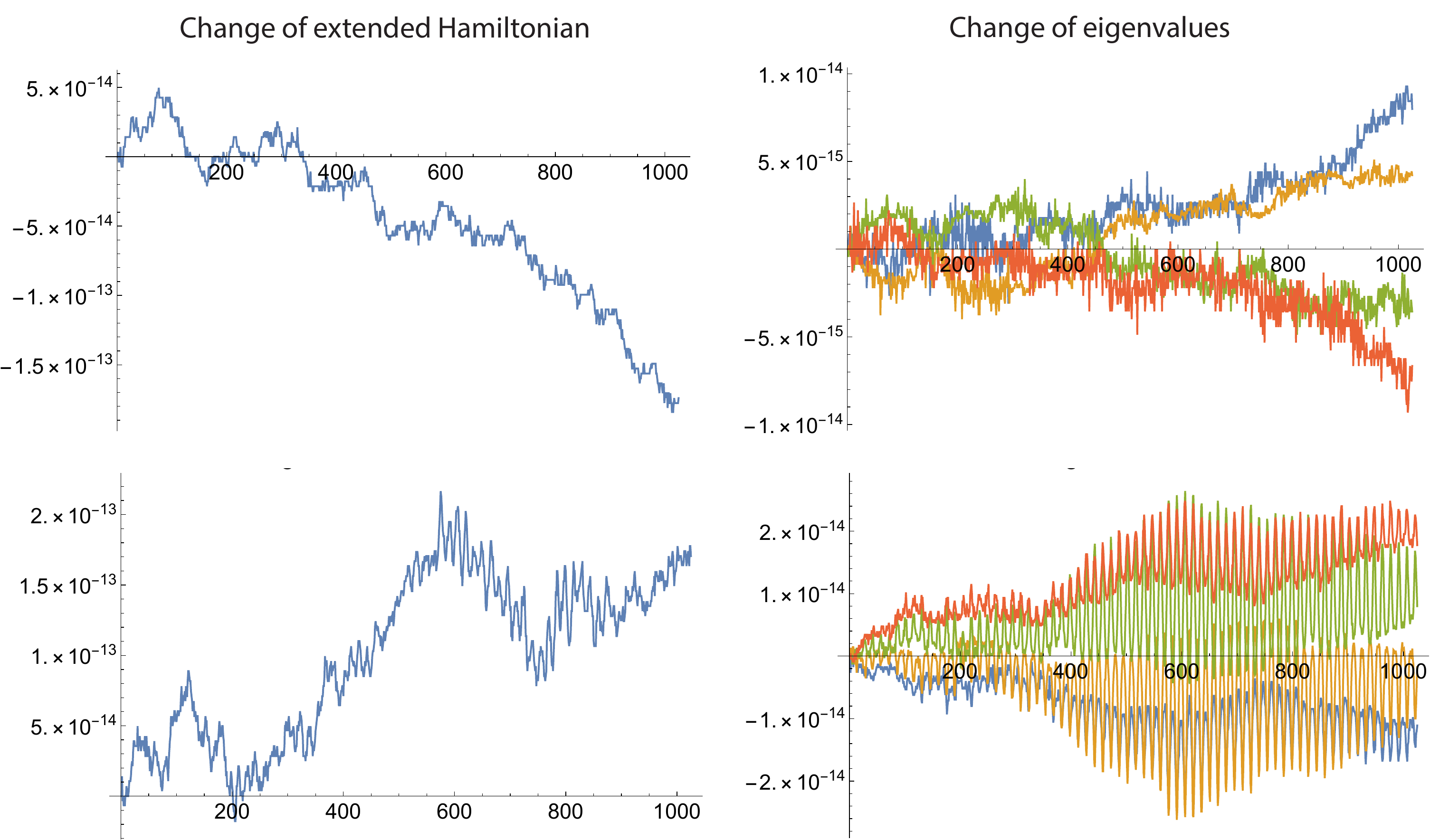}
  \caption{Change in energy (left) and change in eigenvalues (right) for the isospectral implicit midpoint rule and the $4^{\textrm{th}}$-order SDIRK for the Toda lattice.}
  \label{fig:toda_experiments}
\end{figure}
\begin{align*}
  L=
  \begin{pmatrix}
  a_1 & b_1 & 0 & \ldots & b_n \\
  b_1 & a_2 & b_2 & \ldots & 0 \\
  0 & b_2 & a_3 & \ldots & 0 \\
  \vdots & \vdots & \vdots & \ddots & \vdots \\
  b_n & 0 & 0 & \ldots & a_n 
  \end{pmatrix}
  ,\quad
  B(L) =
  \begin{pmatrix}
  0 & b_1 & 0 & \ldots & -b_n \\
  -b_1 & 0 & b_2 & \ldots & 0 \\
  0 & -b_2 & 0 & \ldots & 0 \\
  \vdots & \vdots & \vdots & \ddots & \vdots \\
  b_n & 0 & 0 & \ldots & 0 
  \end{pmatrix} 
\end{align*}
While the Toda lattice is isospectral, it is not Lie-Poisson. 
It can, however, be generalized to an isospectral flow in $\mathfrak{gl}(n)$ by introducing
\begin{align*}
  B(W) =
  \begin{pmatrix}
    0 & W_{12} & 0 & \ldots & -W_{1n} \\ 
    -W_{21} & 0 & W_{23} & \ldots & 0 \\ 
    0 & 0-W_{32}& 0 & \ldots & 0 \\ 
    \vdots & \vdots & \vdots & \ldots & \vdots \\ 
    W_{n1} & 0 & 0 & \ldots & 0 \\ 
  \end{pmatrix} 
\end{align*}
which corresponds to the extended Hamiltonian $\tilde{H}(W) = -\text{Tr} \left( W^\dagger B(W) \right)+ H(W) $. 
The dynamics of the extended Lie-Poisson system are $\dot{W} = [\dd \tilde{H}(W)^\dagger, W]$ which has the form of Eq.~\ref{eq:lie_poisson_algebra}.

In Fig.~\ref{fig:toda_experiments} we show experimental results for the periodic Toda lattice for the implicit midpoint rule and the $7$-stage $4^{\textrm{th}}$-order Runge-Kutta scheme.
In both cases we obtain excellent conservation of the eigenvalues and the extended Hamiltonian. 
The initial condition $W_0$ is the same as used in~\cite{modin2019lie} where $a_i = b_i = (-1)^i, i = 1, \ldots ,4$ and we used $h = 0.1$ for the time step.

\subsection{Matrix model of barotropic vorticity equation on $S^2$}

The barotropic vorticity equation on the sphere is given by
\begin{align}
  \label{eq:euler_equation}
  \dot{\zeta} = \{ \zeta , \Delta^{-1} \zeta \} ,
\end{align}
with $\zeta$ being vorticity and $\{ \, , \}$ the canonical Poisson bracket on $S^2$.
Eq.~\ref{eq:euler_equation} is a simple model for weather and climate dynamics and plays an important role in various applications.

\begin{figure}
  \includegraphics[width=\textwidth]{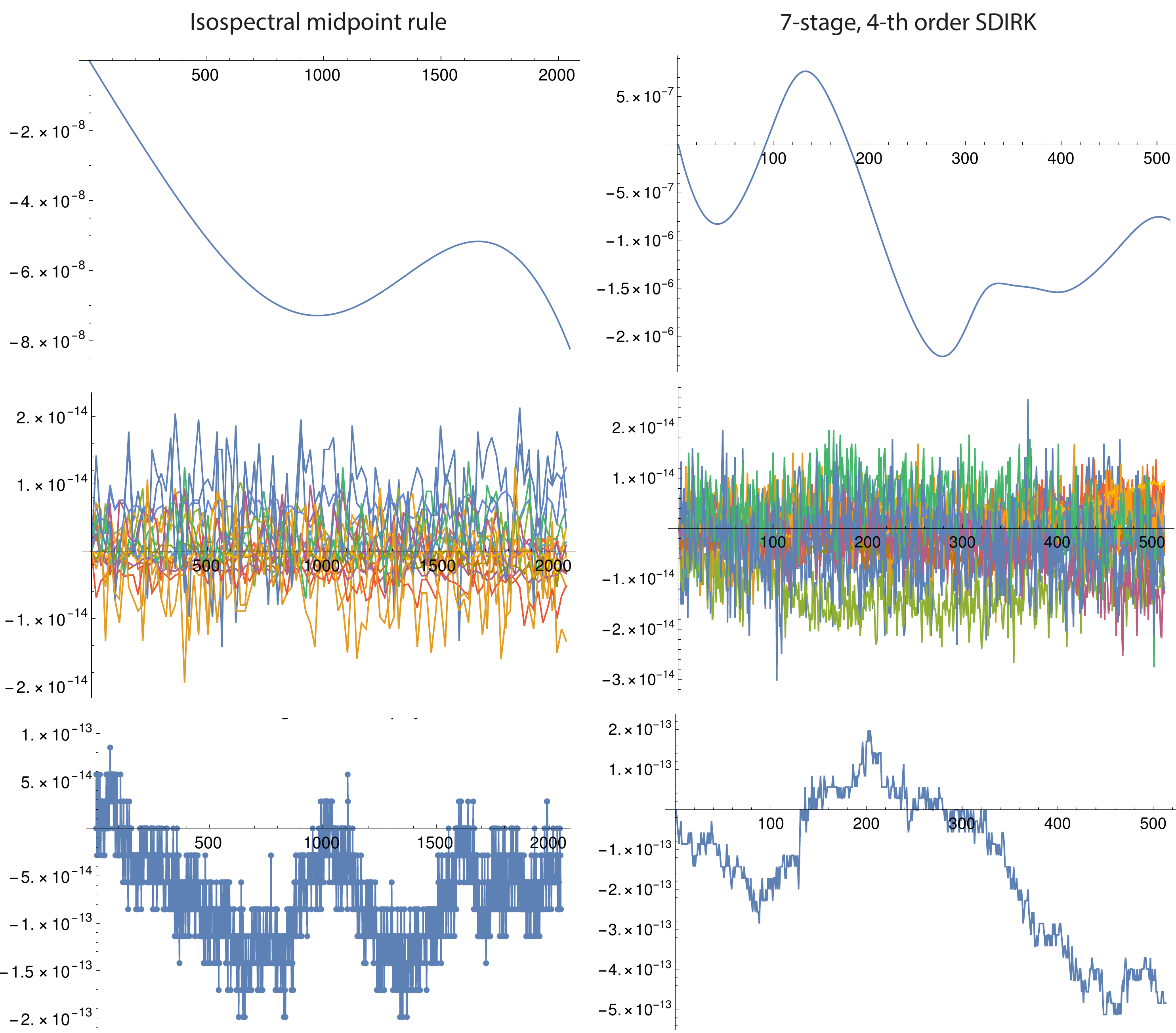}
  \caption{Change in energy (top), eigenvalues (middle) and enstrophy (bottom) for the isospectral implicit midpoint rule and the $4^{\textrm{th}}$-order SDIRK for the matrix model of the barotropic vorticity equation. We used a truncation at $L=16$ ($N = 33$).}
  \label{fig:barotropic_experiments}
\end{figure}

A structure preserving discretization for it is the matrix model by Zeitlin~\cite{Zeitlin2004} that goes back to early work for the torus~\cite{Zeitlin1991} and results by Hoppe on matrix harmonics for $S^2$~\cite{Hoppe1998}.
Recently, Modin and Vivani~\cite{modin2019lie} used it to study turbulence on the sphere.
The governing equation for the matrix model is 
\begin{align}
  \label{eq:barotropic:matrix}
  \dot{W} = N^{3/2} \, [ \hat{\Delta} W , W ]
\end{align}
where $\hat{\Delta}$ is a discrete Laplacian~\cite{Hoppe1998} and $W \in \mathfrak{sl}(N)$  with $N=L+1$ is the representation of vorticity.
It is related to the spatial vorticity field by
\begin{align*}
  W = \sum_{l=1}^L \sum_{m=-l}^l \zeta_{lm} \, T_{lm}^{(L)}
\end{align*}
where the $\zeta_{lm}$ are the spectral coefficients of vorticity with respect to spherical harmonics and the $T_{lm}$ are matrix harmonics~\cite{Hoppe1998}.

Eq.~\ref{eq:barotropic:matrix} is an isospectral flow on $\mathfrak{sl}(N)$ with conserved quantities
\begin{align*}
  C_n = \mathrm{tr}( W^n ) , \quad n = 1 , \ldots, N, 
\end{align*}
which are the discrete analoges to the integrated powers of vorticity that are conserved for the continuous barotropic vorticity equation~\cite{Khesin1998}.

In Fig.~\ref{fig:barotropic_experiments} we report our experimental results for the $2^{\textrm{nd}}$-order implicit midpoint rules and the $4^{\textrm{th}}$-order SDIRK.
The initial condition $W_0$ was a randomly generated matrix in $\mathfrak{sl}(N)$ and we use $h = 0.0025$ (note that the matrix model results in a time rescaling, cf.~\cite[sec. 2.4]{modin2020casimir}).
Both isospectral time integration schemes preserve the eigenvalues of $W$ up to machine precision, as expected, and this carries over to the Casimirs, as exemplified for the discrete enstrophy $C_2$.

\begin{algorithm}[t]
  \label{alg:isospectral_sdirk}
  \caption{Isospectral SDIRK}
  \begin{algorithmic}[1]
    \Statex {\it Input:} $s$, $(b_1, \ldots, b_s), N_{\text{steps}}, h, \mu_0$, 
    \State \vspace*{5.pt} $\mu_{0,r_0} = \mu_0$ 
    
    \For{$n \gets 1, \, N_{\text{steps}} $} \Comment{Time-stepping}
      \For{$i \gets 1, \, s $}
      \State $h_i \gets b_i \, h$
      \State {\bf Solve} $\mu_{n,r_{i-1}} = \left( \text{Id}-\frac{h_i}{2}B\left( \mu_{n,c_i}\right) \right) \mu_{n,c_i} \left(\text{Id}+\frac{h_i}{2} B\left( \mu_{n,c_i}\right) \right)$ for $\mu_{n,c_i}$ \Comment{nonlinear solve} 
      
      \State $\mu_{n,r_i} \gets \left( \text{Id}+\frac{h_i}{2}B\left( \mu_{n,c_i}\right) \right) \mu_{n,c_i} \left(\text{Id}-\frac{h_i}{2} B\left( \mu_{n,c_i}\right) \right)$ \Comment{explicit update}
      \EndFor    
    \EndFor
  \end{algorithmic}
\end{algorithm}

\section{Conclusion}
\label{sec:conclusion}

In this work, we showed that isospectral symplectic diagonally implicit Runge-Kutta methods (SDIRKs) can be derived from the discrete Hamilton's principle of~\cite{BouRabee2007,gawlik2011geometric}.
The resulting numerical schemes are the simplest and most easily implemented isospectral higher order Runge-Kutta methods.
We demonstrate this with our numerical results for the rigid body, the Toda lattice and the discrete 2D Euler fluid.
Our results for a $4^{\textrm{th}}$-order SDIRK yield a time integration that preserve energy and isospectrality almost up to machine precision.
In the available precision, it is hence a faithful discretization of the system.

Some interesting directions for future work are as follow. 
Theorem~\ref{th:variational_generalSDIRK} only applies to SDIRKs and a generalization to arbitrary symplectic Runge-Kutta methods is hence open.
It is currently also not clear to us why the derivation of Gawlik et al.~\cite{gawlik2011geometric} yields the midpoints of the isospectral implicit midpoint rules by Vivani~\cite{viviani2019minimal}. 
We believe that this might shed more light onto the overall structure of discrete variational principles for isospectral systems.
An interesting objective for future work is also to extend our results to systems whose configuration space is a semi-direct product. 
In some cases, such as the heavy top, then also an infinite number of conserved quantities exist. 
Bobenko and Suris~\cite{Bobenko1999} studied this already for the special case of the Lagrange top although integrability is no longer given in the general case.

\begin{acknowledgements}
We would like to acknowledge helpful discussions with Mathieu Desbrun.
\end{acknowledgements}

\appendix
\section{Right invariance}
\label{sec:right_invariance}
The canonical hamiltonian system in Eq.~\ref{eq:hamilton_eq_iso} resulting from the reconstruction of the isospectral system in Eq.~\ref{eq:lie_poisson_algebra} is left invariant. A right invariant version of the canonical hamiltonian system and the algorithm in Eq.~\ref{eq:algorithm} can be achieved as follows.

Consider the left action of the Lie group G on $T^*G$ given by $h\cdot(g,p) = ( hg, (h^\dagger)^{-1})$. It induces the momentum map $\mathcal{J}(g,p) = pg^\dagger$.
Defining the right invariant Hamiltonian on $T^*G$ as $\tilde{H}(g,p) = H \circ \mathcal{J}(g,p)$, where $H$ is the Hamiltonian for the Lie-Poisson system in Eq.~\ref{eq:lie_poisson_algebra}, we have the canonical right invariant Hamiltonian system on $T^*G$
\begin{subequations}
  \begin{align*}
    \dot{g} &= \nabla H(g^\dagger p) g = B(g^\dagger p)^\dagger g \\ \nonumber 
    \dot{p} &= - \nabla H(g^\dagger p)^\dagger p = - B(g^\dagger p) p 
  \end{align*}
\end{subequations}
where again $B(\cdot) = \nabla H (\cdot)^\dagger$.
Applying a symplectic diagonally implicit Runge-Kutta method to the above equation and proceeding as in section \ref{sec:rk_cayley}, we have a right invariant version of the algorithm in Eq.~\ref{eq:algorithm}, 
\begin{subequations}
	\begin{align*}
		\mu_{n, r_{i-1}} &=  \dd \mathrm{cay}^{-1}_{-h_{i}\, B({\mu}_{n,c_i}))}\,\mu_{n,c_i}\\
		\mu_{n, r_{i}} &=  \dd \mathrm{cay}^{-1}_{h_{i}\, B({\mu}_{n,c_i}))}\,\mu_{n,c_i}. 
	\end{align*}
\end{subequations}
Notice that the correct time direction is only obtain when one distinguishes the left and right invariant case.
Theorem~\ref{th:variational_generalSDIRK} also has a right invariant version where the update equation in Eq.~\ref{eq:group_update_generalSDIRK} becomes $g_{n,r_i} = \tau(h_i \xi_{n,c_i}) g_{n,r_{i-1}}$. 		


\printbibliography

\end{document}